\documentclass[11pt]{amsart}
\usepackage{amsmath,amsfonts,amssymb}

\usepackage{latexsym}
\usepackage{eufrak}
\theoremstyle{plain}
\newtheorem{lemma}{Lemma}[section]
\newtheorem{theorem}[lemma]{Theorem}

\theoremstyle{definition}
\newtheorem{remark}[lemma]{Remark}

\numberwithin{equation}{section}

\newcommand{\ZZ}{\mathbb{Z}}

\newcommand{\RR}{\mathbb{R}}

\newcommand{\mA}{\mathcal{A}}

\newcommand{\mT}{\mathcal{T}}

\newcommand{\half}{\frac{1}{2}}

\newcommand{\hP}{\hat{P}}

\newcommand{\ep}{\epsilon}

\begin{document}

\title[Controlled RW with a target site]{Controlled random walk with a target site}
\author{Kenneth S. Alexander}
\address{Department of Mathematics KAP 108 \\
University of Southern California\\
Los Angeles, CA  90089-2532 USA}
\email{alexandr@usc.edu}
\thanks{This research was supported by NSF grant DMS-0804934.}

\keywords{random walk, stochastic control}
\subjclass[2010]{Primary: 60G50; Secondary: 93E20}

\maketitle

\begin{abstract} 
We consider a simple random walk $\{W_i\}$ on $\ZZ^d, d=1,2$, in which the walker may choose to stand still for a limited time.  The time horizon is $n$, the maximum consecutive time steps which can be spent standing still is $m_n$ and the goal is to maximize $P(W_n=0)$.  We show that for $d=1$, if $m_n \gg (\log n)^{2+\gamma}$ for some $\gamma>0$, there is a strategy for each $n$ yielding $P(W_n = 0) \to 1$.  For $d=2$, if $m_n \gg n^\ep$ for some $\ep>0$ then there are strategies yielding $\liminf_n P(W_n=0)>0$.  
\end{abstract}

\section{Introduction.}
We consider a process $\{W_i, 0 \leq i \leq n\}$ on $\ZZ^d$ ($d=1$ or 2) in which $W_0=0$ and each step $W_i - W_{i-1}$ either is 0 (i.e.\ standing still) or is a step $\pm e_i$ of symmetric simple random walk (SSRW), with $e_i$ being the $i$th unit coordinate vector.  The choice between standing still and SSRW step is determined by a strategy.  Formally, a \emph{strategy} (or \emph{n-strategy}) is a mapping $\delta_n: \mT_n \to \{0,1\}$ defined on the space
\[
  \mT_n = \big\{ \{(i,w_i), 0 \leq i \leq j\}: 0 \leq j \leq n-1, |w_i - w_{i-1}| \leq 1 \text{ for all } i \big\}
\]
of all space-time trajectories of length less than $n$; here the values 0 and 1 for $\delta_n$ correspond to standing still and taking a SSRW step, respectively.  Thus the choice of whether to take a step at a time $i$ depends on the trajectory up to time $i-1$.  Letting $\{Y_i, i \geq 1\}$ be SSRW, we then construct the controlled random walk process iteratively from the strategy by
\[
  W_0 = 0, \quad W_i - W_{i-1} = \Delta_i (Y_i - Y_{i-1}) \quad i \leq n,
\]
where 
\[
  \Delta_i = \delta_n\bigg(\{ (k,W_k):  0 \leq k \leq i-1\}\bigg).
\]
We define the time since the last SSRW step to be 
\[
  J_i = \min\{j \in [0,i-1]: W_{i-j} \neq W_{i-j-1}\},
\]
with $J_i = i$ if the set in this definition is empty.  We designate a maximum number $m_n-1$ of consecutive steps standing still and say that a strategy is \emph{admissible} if $J_i \leq m_n-1$ a.s.\ for all $i \leq n$.  A strategy is \emph{Markov} if the value of $\delta_n\big(\{ (k,W_k):  0 \leq k \leq i-1\}\big)$ depends only on the pair $(W_i,J_i)$.  It is easily seen that if the strategy is Markov, then $\{(W_i,J_i): i \geq 0\}$ is a Markov process.  We write $\mA_n$ for the set of all admissible strategies and $\mA_n^*$ for the set of all admissible Markov strategies.  When the dependence on the strategy needs to be made clear, we write $P(\cdot; \delta_n)$ for probability when strategy $\delta_n$ is used, but generally we suppress the $\delta_n$ in the notation.

We are interested in steering the process toward a target site, and specifically in the behavior of
\begin{equation} \label{supdef}
  \hP(W_n=0) = \sup_{\delta_n \in \mA_n} P(W_n=0; \delta_n)
\end{equation}
as $n\to \infty$.  The word ``steering'' is a bit misleading here, as the process never has a drift.  For SSRW, of course $P(W_n=0) \sim Cn^{-d/2}$.  (Here and throughout the paper, $C$ and $C_1,C_2,\dots$ are generic constants, and $a_n \sim b_n$ means the ratio converges to 1.)  A simple strategy to increase this is to minimize steps, i.e.\ always stand still as long as allowed, taking only $n/m_n$ steps, yielding $P(W_n=0) \asymp (m_n/n)^{d/2}$, where we use $a_n \asymp b_n$ to mean the ratio is bounded away from 0 and $\infty$.  A slightly more sophisticated strategy is to minimize steps until time $n-m_n$ then maximize steps (i.e.\ never stand still) until $\{W_i\}$ hits 0 (if it does), and then stand still until time $n$.  For $d=1$, conditionally on $\{|W_{n-m_n}| \leq m_n^{1/2}\}$ this strategy has a success probability of order 1, so the unconditional success probability is of order $\min(m_n/n^{1/2},1)$; in particular it is bounded away from 0 if $m_n \geq Cn^{1/2}$. 

This brings up two questions.  For which $\{m_n\}$ is $\hP(W_n=0)$ bounded away from 0?  And are there nontrivial $\{m_n\}$ for which $\hP(W_n=0) \to 1$?  Our two main theorems give some answers.

\begin{theorem}\label{1d}
Suppose $d=1$ and $m_n \gg (\log n)^{2+\gamma}$ for some $\gamma>0$..  Then $\hP(W_n=0) \to 1$ as $n\to \infty$.
\end{theorem}

In one dimension, the ``more sophisticated'' strategy described above relies on the fact that a SSRW started at distance $k$ from 0 has a probability of order one to hit 0 by time $k^2$.
In two dimensions, this probability is only of order $(\log k)^{-1}$ so it is more difficult to construct a strategy based on waiting for the RW to return to 0 after it has wandered away.  Nonetheless we have the following.

\begin{theorem} \label{2d}
Suppose $d=2$ and $m_n \gg n^{\epsilon}$, for some $\ep>0$.  Then 
\[
  \liminf_n \hP(W_n=0) > 0.
  \]
\end{theorem}

In \eqref{supdef} one could consider only Markov strategies, i.e.\ take the sup over $\mA_n^*$.  Since the underlying SSRW is Markov, one cannot actually do better with a non-Markov strategy, so the two sups are the same.  Allowing non-Markov strategies simply lets us use ones that can be concisely described and analyzed.	

\begin{remark}
One could consider alternate ways of slowing down the controlled RW, in place of standing still.  For example, one could allow a choice between a standard SSRW step and a delayed SSRW step, the latter meaning we stand still with probability $1 - \frac{1}{m_n}$ and take a standard SSRW step of $\pm e_i$ with probability $1/m_n$, with no limit on how many consecutive times we choose the delayed SSRW step.  A brief examination of the proofs shows that both theorem statements above remain valid in this case.
\end{remark}

\begin{remark}
The continuous analog of the problem we consider is a diffusion $\xi_t$ in which the drift is always 0 and one can control the diffusivity $\sigma(x,t)$, but constrained to an interval $[\sigma_1,\sigma_2]$, on a time interval $[0,T]$, where $\sigma_1>0$.  We take $[\sigma_1,\sigma_2] = [\ep_T,1]$ and ask, how slowly can we have $\ep_T \to 0$ as $T \to \infty$ and still have $\hP(|\xi_T|<1) \to 1$ or $\liminf_T P(|\xi_T|<1) >0$?  Note that $\ep_T$ is the analog of $1/m_n$.  McNamara \cite{Mc85} considered closely related questions for a one-dimensional diffusion; see also \cite{Mc83} for another variant.  He proved that there are constants $K=K(\sigma_1,\sigma_2)$ and $\beta = \beta(\sigma_1/\sigma_2)$ such that 
\begin{equation} \label{upperd}
  \hP^x\left(|\xi_T| \leq \frac{\delta}{2} \right) \leq K\left( \frac{\delta}{\sqrt{T}} \right)^\beta \quad \text{for all } x \in\RR, \delta>0,
\end{equation}
while in the other direction, for each $h>0$ there exists $\Delta>0$ such that
\begin{equation} \label{lowerd}
  \hP^0\left(|\xi_T| \leq \frac{\delta}{2} \right) > K\left( \frac{\delta}{\sqrt{T}} \right)^{\beta+h} \quad \text{for all } \delta\in (0,\Delta).
\end{equation}
Here $\hP^x$ denotes probability (maximized over the allowed controls) for a process started at $x$.
Numerical evidence was given that $\beta$ is approximately proportional to $\sigma_1/\sigma_2$.  If we assume this to be true and consider $\ep_T$ of order $1/\log T$, we get $\beta$ of order $1/\log T$ as well.  If we could take $h$ also of this order, then the right side of \eqref{lowerd} would be bounded away from 0 in $T$, as desired.  The problem is that one cannot get a useful result from \eqref{lowerd} if one takes $h$ depending on $T$, since then $\delta$ must also depend on $T$.  Nonetheless we may observe that $\ep_T \leq C/\log T$ corresponds to $m_n \geq C\log n$ in our discrete problem.  Further, still assuming $\beta$ proportional to $\sigma_1/\sigma_2$, we see that since we want the left side of \eqref{upperd} bounded away from 0, we cannot allow $\ep_T \gg 1/\log T$, suggesting that we perhaps cannot do better than requiring $m_n \geq C\log n$ in Theorem \ref{1d}.
\end{remark}

\section{Proof of Theorem \ref{1d}}
Throughout the paper we will make use of various quantities which approach infinity as $n\to \infty$.  
For $\lambda>0$ let $u_n(\lambda) = \max\{k: m_n^{1+k\lambda} \leq n\}$; note that $u_n(\lambda)$ is of order $\log n/\log m_n$ for all $\lambda>0$.  
Also, the hypothesis that $m_n \gg (\log n)^{2+\gamma}$ for some $\gamma>0$ is equivalent to 
\begin{equation} \label{mnhyp1}
  m_n \gg \left( \frac{(\log n)^2}{\log\log n} \right)^{1/(1-\eta)} \quad \text{for some } \eta \in (0,1).
\end{equation}
Fixing such an $\eta$ and writing $u_n$ for $u_n(\eta)$, we observe that \eqref{mnhyp1} is equivalent to
\[
  m_n^{1-\eta}\log m_n \gg (\log n)^2,
\]
hence also to
\[
  m_n^{1+\eta\theta(m_n^{1-\eta}/\log m_n)^{1/2}} \geq n \quad \text{ for all large $n$, for every } \theta>0,
\]
and therefore finally to
\begin{equation} \label{mnhyp}
  u_n = o\left( \left( \frac{ m_n^{1-\eta}}{\log m_n} \right)^{1/2} \right).
\end{equation}
Fix $n$, write $m$ for $m_n$. and define
\[
  \ep_m = \frac{\log m}{m^{1-\eta}}, \quad \text{so} \quad \ep_m = o\left( \frac{1}{u_n^2} \right) \quad \text{as } n \to \infty.
\]
 We define a sequence of windows $\{t_k\} \times I_k, 1 \leq k\leq u_n,$ in space-time, with size decreasing as $k$ increases and the target $(n,0)$ is approached; we then construct a strategy which makes the space-time trajectory of the process pass through all of these windows, with high probability.  Specifically, let 
\[
  t_k = n-m^{1+\eta(u_n-k)}, \quad 2 \leq k \leq u_n, \quad \text{and } t_0 = 0, \quad t_1 = \frac{t_2}{2}, \quad t_{u_n+1} = n,
\]
\[
  N_k = t_k - t_{k-1}, \quad h_{mk} = \left( \ep_m N_{k+1} \right)^{1/2}.
\]
Let $ I_0 = \{0\}$,
\[
  I_k = [-h_{mk},h_{mk}], \quad 1 \leq k \leq u_n, \quad \text{and } I_{u_n+1} = \{0\}.
\]
We want to find a strategy, and choice of $\ep_m$, for which we can show
\begin{equation} \label{claim}
  \max_{k \leq u_n+1} P(W_{t_k} \notin I_k \mid W_{t_{k-1}} \in I_{k-1}) = o\left( \frac{1}{u_n} \right) \quad \text{as } n \to \infty.
\end{equation}
From \eqref{mnhyp} and \eqref{claim} it follows that
\begin{equation} \label{t1}
  P(W_n=0) \geq \prod_{k=1}^{u_n+1} P(W_{t_k} \in I_k \mid W_{t_{k-1}} \in I_{k-1})
    \geq \left( 1 - o\left( \frac{1}{u_n} \right) \right)^{u_n+1} = 1-o(1)
\end{equation}
as $n\to\infty$, proving the theorem.

For \eqref{claim}, fix $1 \leq k \leq u_n+1$ and let $\tau_0 = \tau_0^{(k)} = \min\{t> t_{k-1}: W_t=0\}$. 
Fix $x\in I_{k-1}$ and observe that
\begin{align} \label{options}
  P(W_{t_k} \notin I_k \mid W_{t_{k-1}} = x) &\leq P\left(\tau_0 > t_k \big| W_{t_{k-1}} = x \right) \notag\\
    &\quad + P\left(\tau_0 \leq t_k, |W_{t_k}| > h_{mk} \big| W_{t_{k-1}} = x \right).
\end{align}
To fully specify, and then bound, these probabilities we need to designate a strategy; we do so by describing how the strategy works between $t_{k-1}$ and $t_k$, for general $k$.  For $k \leq u_n$, we begin by taking all SSRW steps (i.e.\ no standing still) from time $t_{k-1}$ until time $\tau_0 \wedge t_k$.  If $\tau_0 > t_k$, we deem the strategy to have failed and we continue in an arbitrary manner, say all SSRW steps.  
If $\tau_0 \leq t_k$, we continue from time $\tau_0$ to $t_k$ by always standing still for the maximum allowed period of time $m$ during the interval $(\tau_0,t_k]$, that is, we take an SSRW step every $m$th time step.  The last standing period is truncated if it would otherwise go beyond time $t_k$.

For $k=u_n+1$, our strategy during $(t_{u_n},t_{u_n+1}] = (n-m,n]$ is  to maximize steps until the the time $\tau_0 = \tau_0^{(u_n+1)}$ when the process first hits 0 (if $\tau_0 \leq n$), then stand still until time $n$.  

We now bound the first term on the right side of \eqref{options}.
From the Reflection Principle we have (for $\ell,h_{mk}$ of opposite even-odd parity):
\begin{equation} \label{reflec}
  P(\tau_0^{(k)}-t_{k-1}>l \mid W_{t_{k-1}}=x ) = P(-|x| < Y_l < |x|),
\end{equation}
and hence for every $\ep>0$,
\begin{align} \label{normalapp}
  P&(\tau_0^{(k)}-t_{k-1}> yx^2 \mid W_{t_{k-1}}=x ) \sim 2\left( \Phi\left( \frac{1}{\sqrt{y}} \right) - \Phi(0) \right) \notag\\
  &\qquad \text{as $|x| \to \infty$, uniformly over } y \in [\ep,\infty).
\end{align}
The left side of \eqref{reflec} is a nondecreasing function of $|x|$, so for all $x \in I_k$, by \eqref{normalapp}
\begin{align} \label{rootbound}
   P\left(\tau_0 > t_k\ \big|\ W_{t_{k-1}} = x \right) &\leq P\left(\tau_0 - t_{k-1} > N_k\ \big|\ W_{t_{k-1}} = h_{mk} \right) \leq C\sqrt{\ep_m}.
\end{align}

Turning to the second term on the right side of \eqref{options}, it is 0 for $k=u_n+1$ so we consider $k \leq u_n$.
We can condition also on $\tau_0$ as follows:  for $t \in (t_{k-1},t_k]$, using Hoeffding's Inequality \cite{Ho63},
\begin{align} \label{cond}
  P\left(|W_{t_k}| > h_{mk}\ \big|\ \tau_0 = t, W_{t_{k-1}} = x \right) &= P\left( |W_{t_k}| > h_{mk}\ \big|\ W_t = 0 \right) \notag\\
  &= P\left( |Y_{(t_k - t)/m}| > h_{mk} \right) \notag\\
  &\leq e^{-h_{mk}^2m/2N_k} \notag\\
  &= e^{-m\ep_mN_{k+1}/2N_k}.
\end{align}
For $k\geq 1$ we have $N_{k+1}/N_k \geq m^{-\eta}$ (with equality for $k\geq 3$), so \eqref{cond} says that
\begin{equation} \label{cond2}
  P\left(|W_{t_k}| > h_{mk}\ \big|\ \tau_0 = t, W_{t_{k-1}} = x \right) \leq m^{-1/2} \leq \sqrt{\ep_m}.
\end{equation}   

Since $x \in I_{k-1}$ is arbitrary, combining \eqref{options}, \eqref{rootbound} and \eqref{cond2} yields that for all $1 \leq k \leq u_n+1$,
\begin{equation} \label{}
  P(W_{t_k} \notin I_k \mid W_{t_{k-1}} \in I_{k-1}) \leq C\sqrt{\ep_m} = o\left( \frac{1}{u_n} \right) \quad \text{as } n \to \infty,
\end{equation}
proving \eqref{claim} and thereby proving Theorem \ref{1d}.

\section{Proof of Theorem \ref{2d}}
We keep the same definition of $u_n(\lambda)$ and note that now our hypothesis on $m_n$ is equivalent to the statement that $\{u_n\}$ is bounded, say $u_n \leq u< \infty$ for all $n$.  We keep the same formula for $t_k$ but with $\eta$ replaced by $\kappa$, determined as follows.  Choose
\[
  \theta \in \left( \frac{1-\ep}{2},\half \right),
\]
so $1-2\theta < \ep$, then choose $\kappa$ small enough so
\begin{equation} \label{kappa}
  0 < \frac{1-2\theta}{1 - 2\kappa\theta} < \ep.
\end{equation}
We then write $u_n$ for $u_n(\kappa)$.
For our windows, in place of the interval $I_k$ we have the square $Q_k = [-N_{k+1}^\theta,N_{k+1}^\theta]^2$.  
To distinguish dimensions clearly, we now write $\{Y_i^{(d)}\}$ for $d$-dimensional SSRW, $d=1,2$.
We use the same strategy as in one dimension:  in each interval $(t_{k-1},t_k]$, take an SSRW step every time step until time $\tau_0^{(k)}$, then an SSRW step every $m$th time step from time $\tau_0^{(k)}$ to $t_k$.  In place of \eqref{claim}, we will need that for some $C>0$,
\begin{equation} \label{claim2}
  P(W_{t_k} \in Q_k \mid W_{t_{k-1}} \in Q_{k-1}) \geq C \quad \text{for all $k$, for all large } n.
\end{equation}
Let $L_{[a,b]}$ be the number of visits to 0 by the SSRW $\{Y_i^{(2)}\}$ during the time interval $[a,b]$.  
In comparison to \eqref{rootbound}, we have for all $x \in Q_{k-1}$:
\begin{align} \label{rootbound2}
  P\left(\tau_0 \leq t_k\ \big|\ W_{t_{k-1}} = x \right) &= P\left(\tau_0 \leq t_k\ \big|\ Y_{t_{k-1}}^{(2)} = x \right) \notag\\
  &\geq \frac{ E(L_{(t_{k-1},t_k]} \mid Y_{t_{k-1}}^{(2)} = x) }{ E(L_{(t_{k-1},t_k]} \mid \tau_0 \leq t_k, Y_{t_{k-1}}^{(2)} = x) } 
\end{align}
Since $|x| \leq CN_k^\theta$ with $\theta<\half$, the numerator in \eqref{rootbound2} is equal to
\begin{align} \label{num}
  \sum_{i=1}^{N_k} P(Y_i=0 \mid Y_0^{(2)} = x) \geq C\sum_{i=|x|^2}^{N_k} \frac{1}{i} \geq C_1\log N_k.
\end{align}
Here $C_1$ depends on $\theta$.
The denominator in \eqref{rootbound2} is bounded above by
\begin{align} \label{denom}
  E(L_{(t_{k-1},t_k]} \mid Y_0^{(2)}=0) \leq C\sum_{i=1}^{N_k} \frac{1}{i} \leq C_2\log N_k.
\end{align}
Therefore we have the analog of \eqref{rootbound}:
\begin{equation} \label{rootbound3}
  P\left(\tau_0 \leq t_k\ \big|\ W_{t_{k-1}} = x \right) \geq \frac{C_1}{C_2} >0 \quad \text{for all } k \leq u_n, x \in Q_{k-1}.
\end{equation}

%Let us write the SSRW as $Y_j = (Y_j^{(1)},Y_j^{(2)})$.  
In comparison to \eqref{cond}, we have for $x \in Q_k$ and $t \in (t_{k-1},t_k]$, using again Hoeffding's Inequality \cite{Ho63}:
\begin{align} \label{cond3}
  P\left(W_{t_k} \notin Q_k\ \big|\ \tau_0 = t, W_{t_{k-1}} = x \right) &= P\left( W_{t_k} \notin Q_k\ \big|\ W_t = 0 \right) \notag\\
  &= P\left( Y_{(t_k - t)/m} \notin Q_k \mid Y_0^{(2)}=0 \right) \notag\\
  &\leq 2\max_{j \leq N_k/m} P(|Y_j^{(1)}| > N_{k+1}^\theta) \notag\\
  &\leq 4e^{-mN_{k+1}^{2\theta}/2N_k}.
\end{align}
For $k \geq 1$ we have $N_{k+1} \geq m^{-\kappa}N_k$ (with equality for $k \geq 3$), so by \eqref{kappa},
\begin{equation} \label{exponent}
  \frac{mN_{k+1}^{2\theta}}{N_k} \geq \frac{m^{1-2\kappa\theta}}{N_k^{1-2\theta}} \geq \frac{m^{1-2\kappa\theta}}{n^{1-2\theta}}
    \to \infty \text{ as } n \to \infty.
\end{equation}
Combining \eqref{rootbound3}, \eqref{cond3} and \eqref{exponent} we see that \eqref{claim2} holds.  Since $\{u_n\}$ is bounded, it follows as in \eqref{t1} that $\liminf_n P(W_n=0) > 0$.

\section{Acknowledgements}

The author would like to thank Ananda Weerasinghe for helpful conversations.

\end{document}